\documentclass[a4paper,12pt]{article}
\usepackage[english]{babel}
\usepackage{mathrsfs}
\usepackage{epsfig,afterpage}
\usepackage{psfrag}
\usepackage{graphicx}
\usepackage{amsmath}
\usepackage{amssymb}
\usepackage{amsfonts}
\usepackage{multirow}

\newtheorem{theorem}{Theorem}

\newtheorem{remark}{Remark}

\newtheorem{exam}{Example}

\def\{{\protect\lbrace}
\def\}{\protect\rbrace}
\def\Real{\operatorname{Re}}
\def\dim{\operatorname{dim}}
\def\const{{\rm{const}}}
\def\div{\operatorname{div}}

\def\alf{\alpha}
\def\bet{\beta}
\def\Del{\Delta}
\def\eps{\varepsilon}
\def\gam{\gamma}
\def\Gam{\Gamma}
\def\Gamo{\Gam_0}
\def\phi{\varphi}
\def\om{\omega}

\def\si{\sigma}
\def\la{\lambda}
\def\ov{\overline}

\let\<\langle
\let\>\rangle

\newcommand{\bR}{\mathbb{R}}

\newcommand{\bZ}{\mathbb{Z}}

\newcommand{\ox}{\bar x}
\newcommand{\oy}{\bar y}
\newcommand{\op}{\bar p}
\newcommand{\FF}{\mathscr {F}}

\newcommand{\gama}{\gam_{\alf}}
\newcommand{\gamap}{\gam_{\alf}^+}
\newcommand{\gamam}{\gam_{\alf}^-}
\newcommand{\gamapm}{\gam_{\alf}^{\pm}}

\newcommand{\gamabp}{\gam_{\alf, \bet}^+}

\newcommand{\PR}{{\mathbb R}{\mathbb P}}
\newcommand{\pa}{\partial}

\begin{document}

\title{A brief survey on singularities of geodesic flows \\ in smooth signature changing metrics on 2-surfaces}

\author{
N.G. Pavlova\footnote{
Department of Nonlinear Analysis and Optimization, 
RUDN University (Moscow, Russia).
Email: natasharussia@mail.ru.
}
\, and A.O. Remizov\footnote{
CMAP \'{E}cole Polytechnique (Palaiseau, France), CNRS.
Email: alexey-remizov@yandex.ru.
}
}

\maketitle

\begin{abstract}
We present a survey on generic singularities of geodesic flows in smooth signature changing metrics
(often called {pseudo-Riemannian}) in dimension~2.
Generically, a pseudo-Riemannian metric on a 2-manifold $S$
changes its signature (degenerates) along a curve $S_0$, which locally separates
$S$ into a Riemannian ($R$) and a Lorentzian ($L$) domain.
The geodesic flow does not have singularities over $R$ and $L$, and
for any point $q \in R \cup L$ and every tangential direction $p \in \PR$ there exists a unique geodesic
passing through the point $q$ with the direction $p$.
On the contrary, geodesics cannot pass through a point
$q \in S_0$ in arbitrary tangential directions, but only in some admissible directions;
the number of admissible directions is 1 or 2 or 3.
We study this phenomenon and the local properties of geodesics near $q \in S_0$.
\end{abstract}

{
\renewcommand{\thefootnote}{\fnsymbol{footnote}}
\footnote[0]{2010 Mathematics Subject classification 53C22, 53B30, 34C05.} 
}

\section{Introduction}

Let $S$ be a real smooth manifold, $\dim S = n \ge 2$.
By \textit{metric} on $S$ we mean a symmetrical covariant tensor field of the second order on the tangent bundle $TS$, not necessary positive defined. Moreover, metrics whose signature has different signs at different points of $S$, are of the special interest.
For instance, in the quantum theory of gravitation and general relativity two types
of signature changing metrics are considered:
\begin{itemize}
\item
{\bf Smooth}.
The metric is degenerate on a hypersurface $S_0 \subset S$ that divides the \textit{Riemannian} region
$R \subset S$ with signature $(+ \cdots + +)$ from the \textit{Lorentzian} region $L \subset S$ with signature $(+ \cdots + -)$. Example: $ds^2 = dx_1^2 + \cdots + dx_{n-1}^2 + x_n dx_n^2$.
\item
{\bf Discontinuous}.
The metric is smooth and non-degenerate everywhere except for a hypersurface $S_0 \subset S$ (which separates $R$ and $L$ defined as above), where it fails to be continuous.
Example: $ds^2 = dx_1^2 + \cdots + dx_{n-1}^2 + \frac{1}{x_n} dx_n^2$.
\end{itemize}

\begin{figure}[h]
\begin{center}
\includegraphics[width=225pt,height=165pt]{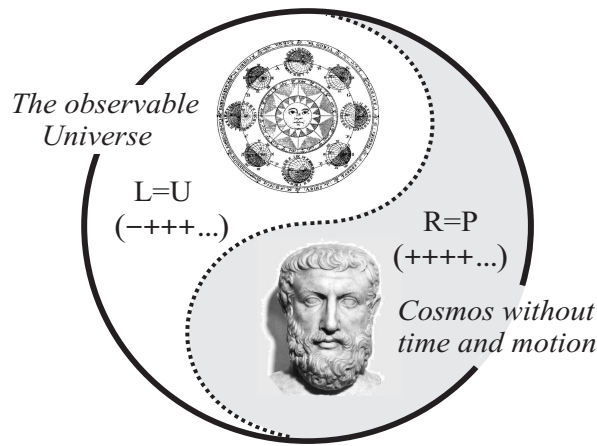}
\caption{
Sakharov's cosmological model. Here the region $L$ (which apparently includes the observable Universe) is denoted by $U$ and the region $R$ (cosmos without time and motion) is denoted by $P$ (after Parmenides, a Greek philosopher theorized about space and time).
The hypersurface $S_0$ is represented by the dotted line.
}
\label{fig1}
\end{center}
\end{figure}

In the paper \cite{Sakhar}, Russian physicist A.D.~Sakharov conjectured there exist states of the physical continuum which include regions with different signatures of the metric; the observed Universe and an infinite number of other Universes arose as a result of quantum transitions with a change in the signature of the metric.
This concept is exemplified by Fig.~\ref{fig1}.

In his cosmological model, Sakharov used discontinuous metrics.
However, some other authors consider models with smooth signature changing metrics;
see e.g. \cite{AFL, KosKri, KosKri1, KosKri2} and the references therein.
From physical viewpoint, the difference between smooth and discontinuous signature changing metrics corresponds to different physical proposals, in particular, different solutions of the Einstein equation.
{\it Euclidean-Lorentzian transitions} (junctions) between the domains $R$ and $L$ play
an important role, both in the smooth and discontinuous models.
The term {\it Euclidean} is used in sense of {\it Riemannian}, that is typical for physical literature,
see e.g. \cite{AB}. Similarly, the term {\it Lorentzian} is referred to non-degenerate indefinite metrics.

In this paper, we discuss a purely mathematical problem connected with smooth signature changing metrics
(further called {\it pseudo-Riemannian}): the local behavior of geodesics in a neighborhood of the points
where the metric has a generic degeneracy.
Such points are singular points of the geodesic flow, and the standard existence and uniqueness theorem
for ordinary differential equations is not applicable.
This leads to an interesting geometric phenomenon: geodesics cannot pass through a degenerate point
in arbitrary tangential directions, but only in certain directions said to be {\it admissible}.

A study of this phenomenon for two-dimensional pseudo-Riemannian metrics is started in
\cite{GR, Rem-Pseudo, Rem15, Rem-Tari};
similar results in three-dimensional case were announced in \cite{PR11}.
In these works, mainly the local properties of geodesics and geodesic flows were considered,
some global properties of geodesics of pseudo-Riemannian metrics
with differentiable groups of symmetries are investigated in \cite{Rem15}.
This allows, in particular, to obtain the phase portraits of geodesics on surfaces of revolution
(sphere, torus, etc.) embedded in three-dimensional Minkowski space.

Various other aspects of pseudo-Riemannian metrics (including the Gauss--Bonnet formula)
are treated by many authors, see e.g. \cite{GKT, KT, KosKri1, KosKri2, Kossowski, Mier, Steller}
and the references therein.
However, there exist a number of unsolved problem connected with degeneracy of metrics.
According to our knowledge, the problem of local geodesic equivalence of pseudo-Riemannian metrics
at degenerate points is not studied yet, although it is well studied for Riemannian and Lorentzian metrics,
see e.g. \cite{Bolsinov} (in this paper, the authors call {\it pseudo-Riemannian} what we call
{\it Lorentzian}, i.e., non-degenerate indefinite metrics).

\medskip

From now we always assume that $\dim S = 2$.

\begin{figure}[h]
\begin{center}
\includegraphics[width=420pt,height=140pt]{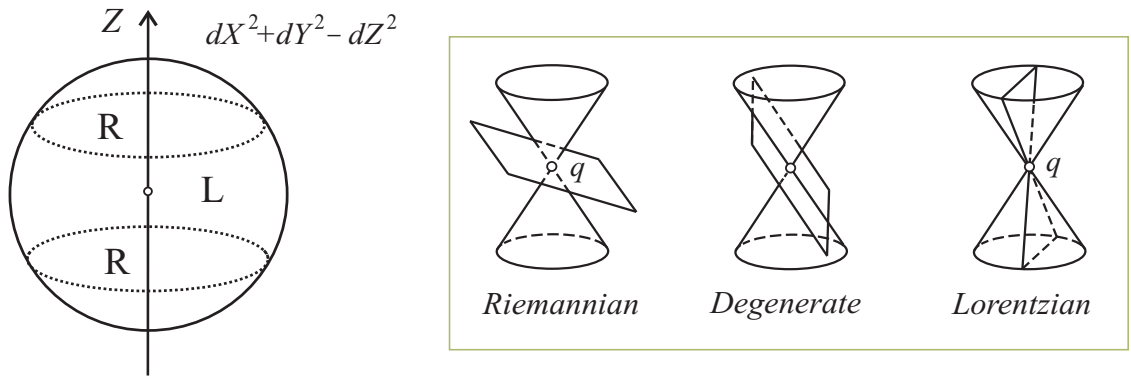}
\caption{
On the left: pseudo-Riemannian metric on the sphere $X^2 + Y^2 + Z^2 = 1$ in 3D Minkowski space.
Here $S_0$ consists of two parallels $Z = \pm 1/{\sqrt{2}}$ depicted as dotted lines.
On the right: intersections of the light cone with the tangent plane $T_q S$, $q \in S$.
}
\label{fig2}
\end{center}
\end{figure}

Similarly, just as Riemannian metrics naturally appear on surfaces embedded in Euclidean space, pseudo-Riemannian metrics can be generated in pseudo-Euclidean space.
Let $S$ be a smooth surface embedded in 3D Minkowski space $(X,Y,Z)$ with the
pseudo-Euclidean metric $dX^2 + dY^2 - dZ^2$.
Then the pseudo-Euclidean metric in the ambient $(X,Y,Z)$-space induced a pseudo-Riemannian metric on $S$.
For instance, let $S$ be the standard Euclidean sphere
$$
X^2+Y^2+Z^2=1.
$$
The metric induced on the sphere $S$ degenerates on two parallels $Z = \pm 1/{\sqrt{2}}$, which
separate $S$ into three regions, where the metric has constant signatures.
The North $\bigl( Z > 1/{\sqrt{2}}\bigr )$ and the South $\bigl( Z < -1/{\sqrt{2}} \bigr )$ regions are
Riemannian, while the equatorial region $|Z|< 1/{\sqrt{2}}$ is Lorentzian; see Fig.~\ref{fig2} (left).
The condition of the point $q \in S$ belonging to $R$ or $S_0$ or $L$ depends on the
mutual relationships between the tangent plane $T_q S$ and the {\it isotropic} ({\it light}) {\it cone}
$$
dX^2+dY^2-dZ^2 = 0;
$$
see Fig.~\ref{fig2} (right).

\eject

\section{Definition of geodesics}

Consider a two-dimensional manifold (surface) $S$ with pseudo-Riemannian metric
\begin{equation}
ds^2 = a(x,y) \, dx^2 + 2b(x,y) \, dx dy + c(x,y) \, dy^2,
\label{1}
\end{equation}
whose coefficients are smooth (i.e., $C^{\infty}$).
Geodesics in the metric \eqref{1} can be defined via variational principles similarly to the Riemannian case,
with additional nuances.

For instance, the arc-length parametrization is not defined for the {\it isotropic lines}
(or {\it lightlike lines} or {\it null curves}).
Moreover, the Lagrangian of the length functional
$$
J_l(\gamma) = \int\limits_{\gamma} \sqrt{a{\dot x}^2 + 2b{\dot x}\dot y + c{\dot y}^2}\,dt \ \to \ {\rm extr},
$$
where the dot means differentiation by the parameter $t$,
fails to be differentiable on the {\it isotropic surface} $\FF$
\begin{equation}
a(x,y) \, dx^2 + 2b(x,y) \, dx dy + c(x,y) \, dy^2 = 0,
\label{2}
\end{equation}
and the Euler--Lagrangian equation for the length functional is not defined on $\FF$.
Note that equation \eqref{2} defines the isotropic surface $\FF$ in the complement of the zero section of $TS$
or, equivalently, in the projectivized tangent bundle $PTS$.

{\it Binary} differential equation \eqref{2} defines a direction field on $\FF$,
whose integral curves correspond to isotropic lines in the metric \eqref{1}.
This equation plays an important role for understanding the behavior of geodesics,
and we consider it in more detail below.

\medskip

As already mentioned above, the Euler--Lagrangian equation for the length functional $J_l$ does not allow to define extremals on $\FF$. However, this problem does not arise if we define geodesics as extremals of the action functional
$$
J_a(\gamma) =
\int\limits_{\gamma} (a{\dot x}^2 + 2b{\dot x}\dot y + c{\dot y}^2) \,dt \ \to \ {\rm extr}.
$$
The corresponding Euler--Lagrange reads
\begin{equation}
\left \{ \
\begin{aligned}
& 2(a \ddot x + b \ddot y) = (c_x-2b_y) {\dot y}^2 - 2a_y {\dot x} {\dot y} -a_x {\dot x}^2, \\
& 2(b \ddot x + c \ddot y) = (a_y-2b_x) {\dot x}^2 - 2c_x {\dot x} {\dot y} - c_y {\dot y}^2, \\
\end{aligned}
\right.
\label{ELA}
\end{equation}
and the corresponding parametrization is called {\it natural} or {\it canonical}.
Obviously, the definition of geodesics as auto-parallel curves in the Levi--Civita connection generated
by the metric \eqref{1} leads to the same Equation \eqref{ELA}.

\begin{figure}[h]
\begin{center}
\includegraphics[width=200pt,height=140pt]{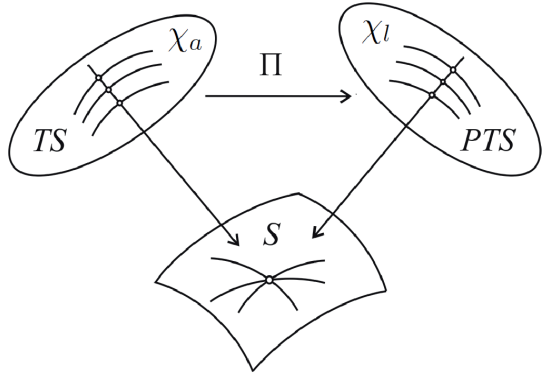}
\caption{
The relationship between the fields $\chi_l$ and $\chi_a$.
The horizontal arrow means the projectivization $\Pi \colon TS \to PTS$.
Two lower arrows are the standard projections $TS \to S$ and $PTS \to S$.
}
\label{fig3}
\end{center}
\end{figure}

The natural parametrization is well defined for all types of geodesics, including isotropic.
For non-isotropic geodesics it coincides with the arc-length (of course, here the length to be real or imaginary).
The functionals $J_l$ (length) and $J_a$ (action) define the corresponding fields of extremals:
$\chi_l$ on $PTS$ away of $\FF$ and $\chi_a$ on the complement of the zero section of $TS$
(including~$\FF$).
The relationship between the fields $\chi_l$ and $\chi_a$ is as follows (see also Fig.~\ref{fig3}).

The natural projectivization $\Pi \colon TS \to PTS$ sends the field $\chi_a$ to a direction field on $PTS$,
which is parallel to the vector field
\begin{equation}
\vec V = 2\Del \biggl(\frac{\pa}{\pa x} + p \frac{\pa}{\pa y} \biggr) + M \frac{\pa}{\pa p},
\ \quad  p = \frac{dy}{dx},
\label{3}
\end{equation}
where
$$
\Del (x,y) = ac-b^2, \ \ \
M(x,y,p) = \sum\limits_{i=0}^{3} \mu_i(x,y)p^i,
$$
with the coefficients
\begin{equation}
\begin{aligned}
&\mu_0 = a(a_y-2b_x) +a_xb,\\
&\mu_1 = b(3a_y-2b_x) + a_xc - 2ac_x, \\
&\mu_2 = b(2b_y-3c_x) + 2a_yc - ac_y, \\
&\mu_3 = c(2b_y-c_x) - bc_y.
\end{aligned}
\label{mu}
\end{equation}

The vector field $\vec V$ given by \eqref{3} is defined and smooth at all points of $PTS$
including the isotropic surface $\FF$.
It is worth observing that the direction field $\chi_l$ is parallel to \eqref{3} at all points where $\chi_l$
is defined, i.e., at all points away from the surface $\FF$.
One can interpret the direction field given by \eqref{3} as a natural extension of $\chi_l$ to $\FF$.
This brings us to the following definition:
the projections of integral curves of the field \eqref{3} from $PTS$ to $S$
distinguished from a point are {\it non-parametrized geodesics}
in the pseudo-Riemannian metric \eqref{1}.

Moreover, let $\vec W$ be the vector field on $PTS$
(determined uniquely up to multiplication by a non-vanishing scalar factor)
that corresponds to the length functional $J_l$.
Since the length functional is invariant with respect to reparametrizations, one can put $t=x$ and
take as $\vec W$ the vector field corresponding to the Euler--Lagrange equation with the Lagrangian
$\sqrt{F}$, where $F(x,y,p) = a(x,y) + 2b(x,y) p + c(x,y) p^2$.
A straightforward calculation (see \cite{Rem15}) shows that
\begin{equation}
\vec W = \frac{1}{2F^{\frac{3}{2}}} \vec V \quad \textrm{and} \quad
\div \vec W = 0 \ \ \textrm{at all points where} \ \ F \neq 0.
\label{5-Oct-2016}
\end{equation}

The field $\vec W$ is divergence-free, since it comes directly from an Euler--Lagrange equation, while
$\vec V$ is not, since it is obtained via an additional procedure, the projectivization $\Pi \colon TS \to PTS$.
The property \eqref{5-Oct-2016} plays an important role, due to the following general fact:

\begin{theorem}[\cite{GR}]
\label{T0}
Let $\vec V(\xi)$, $\xi \in \bR^n$, be a smooth vector field, $f(\xi)$ be a smooth scalar function such that
the hypersurface $\FF = \{ \xi: f(\xi)=0 \}$ is regular, $r$ be a positive real number.
Suppose that the field $\vec W (\xi) = f^{-r}(\xi) \vec V(\xi)$ is divergence-free at all points where it is defined,
i.e.,  at all points $\xi \notin \FF$.
Then $\FF$ is an invariant hypersurface of the field $\vec V$.
Moreover, let $\xi_* \in \FF$ be a singular point of $\vec V$ and $\lambda_1, \ldots, \lambda_n$
be the eigenvalues of the linearization of $\vec V$ at $\xi_*$.
Then $\lambda_1 + \cdots + \lambda_n = r\lambda_j$ for at least one~$j$.
\end{theorem}

By Theorem~\ref{T0}, we have the following assertions:
\begin{itemize}
\item
The isotropic surface $\FF$ is an invariant surface of the field \eqref{3} and all isotropic lines are geodesics
(with identically zero length).\footnote{
The first assertion is valid for any $\dim S \ge 2$,
while the second assertion (about isotropic lines) is valid for $\dim S = 2$ only.
Indeed, in the case $\dim S > 2$ there exist isotropic lines that are not geodesics;
see the example in \cite{Rem15}.
}
\item
Geodesics do not change their type (timelike, spacelike, isotropic) away of degenerate points.
This statement follows from the previous one.
\end{itemize}


\section{Equation of isotropic lines}

Suppose that set
$$
S_0 = \{ q=(x,y) \in S \colon \Del(x,y) = 0 \}
$$
is a regular curve. It is called the {\it degenerate} or {\it discriminant} curve of the metric \eqref{1},
and points $q \in S_0$ are called {\it degenerate points} of the metric.
Then the coefficients $a, b, c$ do not vanish simultaneously, and the isotropic direction
\begin{equation}
p_0(q) = -\frac{a}{b}(q) = -\frac{b}{c}(q), \quad q \in S_0,
\label{40}
\end{equation}
is defined and unique at every point $q \in S_0$.

\begin{figure}[h]
\begin{center}
\includegraphics[height=5.9cm]{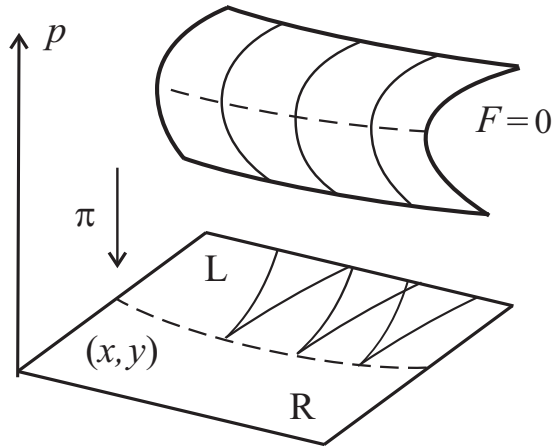}
\caption{
The isotropic surface $\FF$ in $PTS$, integral curves of the field $X$ (top) and integral curves of the
equation $F(x,y,p)=0$ (down).
The dashed lines represent the criminant (top) and the discriminant curve (down).
}
\label{fig0}
\end{center}
\end{figure}

The projectivization $\Pi \colon TS \to PTS$ transforms binary differential equation \eqref{2}
into the implicit differential equation
\begin{equation}
F(x,y,p)=0,  \quad \textrm{where} \quad F = a(x,y) + 2b(x,y) p + c(x,y) p^2.
\label{IDE}
\end{equation}
In the space $PTS$, the surface $\FF$ composes a two-sheeted covering of
the Lorentzian domain of $S$ ($\Del<0$) with branching along the discriminant curve $S_0$.
Over the Riemannian domain ($\Del>0$), the surface $\FF$ does not pass. See Fig.~\ref{fig0}.

A well-known geometrical approach to study implicit equation \eqref{IDE} consists of the lift
the multivalued direction field on $S$ to a single-valued direction field $X$ on the surface $\FF$.\footnote{
This approach is applicable to implicit differential equations $F(x,y,p)=0$ with a smooth function $F$
not necessarily quadratic in $p$. The idea goes back to H.\,Poincar\'e and A.\,Clebsch, see \cite{Rem06} for details.
}
The field $X$ is an intersection of the contact planes $dy = pdx$ with
the tangent planes to the surface $\FF$, that is, $X$ is defined by the vector field
\begin{equation}
\dot x = F_p, \ \ \  \dot y = pF_p, \ \ \ \dot p = -(F_x+pF_y),
\label{6-Oct-2016}
\end{equation}
whose integral curves become isotropic lines of the metric \eqref{1} after the projection $\pi: \FF \to S$
along the $p$-direction. Further we shall call this direction {\it vertical} in the space $PTS$.
The locus of the projection $\pi: \FF \to S$ (given by the equations $F=F_p=0$) is called
the {\it criminant} of equation \eqref{IDE}. It is not hard to see that the criminant
consists of the points $(q,p_0(q))$, $q \in S_0$ (see formula \eqref{40}).

\begin{figure}[h]
\begin{center}
\includegraphics[height=130pt]{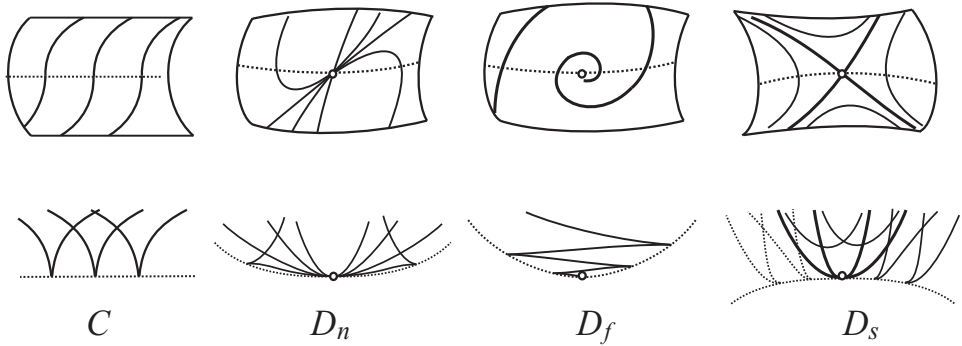}
\caption{
From the top to the bottom:
integral curves of the field \eqref{6-Oct-2016} on the isotropic surface $\FF$ and
isotropic lines, obtained by the the projection $\pi: \FF \to S$.
The dashed lines represent the criminant (top) and the discriminant curve (down).
}
\label{fig5}
\end{center}
\end{figure}

Since $\FF$ is an invariant surface of the field \eqref{3} and both fields \eqref{3} and \eqref{6-Oct-2016} are
tangent to the contact planes $dy = pdx$, the restriction of \eqref{3} to the invariant surface $\FF$
is parallel to \eqref{6-Oct-2016}.
Moreover, the restriction of the field \eqref{3} to $\FF$ is equal to the field \eqref{6-Oct-2016} multiplied by a smooth scalar function vanishing along the criminant (see \cite{GR}).
Generically, here there are two possible cases:
\begin{itemize}
\item
The case $C$: the isotropic direction $p_0(q)$ is transversal to $S_0$.
Then the field \eqref{6-Oct-2016} at the point $(q,p_0(q))$, $q \in S_0$, is non-singular,
and binary equation \eqref{2} has {\it Cibrario normal form} $dx^2 = y\, dy^2$.
See Fig.~\ref{fig5} (left).
\item
The case $D$: the isotropic direction $p_0(q)$ is tangent to $S_0$.
The field \eqref{6-Oct-2016} at $(q,p_0(q))$, $q \in S_0$, has a non-degenerate singular point:
saddle or node or focus (subcases $D_s, D_n, D_f$, respectively).
Under certain additional conditions (formulated below),
binary equation \eqref{2} has {\it Dara--Davydov normal form}
\begin{equation}
dy^2 = (y-\eps x^2) \, dx^2,
\label{DD}
\end{equation}
where $\eps < 0$ (if saddle) or $0< \eps < \tfrac{1}{16}$ (if node) or $\eps > \tfrac{1}{16}$ (if focus).
See Fig.~\ref{fig5}.
\end{itemize}

The normal form $dx^2 = y\, dy^2$ is named after Italian mathematician Maria Cibrario
who established it first in $C^{\om}$ (real analytic) category when studying second-order linear
partial differential equations of the mixed type \cite{Cib}.
Later on, a general (and rather simple) proof of the Cibrario normal form (in $C^{\om}$ and $C^{\infty}$ categories)
was presented in the famous Arnold's book \cite{A-Geom}.

The normal form \eqref{DD} was firstly conjectured by Brazilian mathematician Lak Dara \cite{Dara}
and then proved by A.A.~Davydov \cite{Dav85} under the following genericity conditions.
Let $\alf_{1,2}$ be the eigenvalues of the linearization of the vector field \eqref{6-Oct-2016}
at the singular point considered.
Then $\alf_{1,2}$ are roots of the characteristic equation $\alf^2 - \alf + 4\eps = 0$,
and the excluded values $\eps =0$ and $\eps = \tfrac{1}{16}$ correspond to a degenerate singular point
(saddle-node or degenerate node, respectively).
The additional conditions required for the normal form \eqref{DD} are the following.

First, the ratio of $\alf_{1,2}$ is different from $\pm 1$,
and the eigendirections are not tangent to the criminant.
Second, the germ of the vector field \eqref{6-Oct-2016} is $C^{\infty}$-linearizable, i.e., it is
$C^{\infty}$-smoothly equivalent to its linear part.
The $C^{\infty}$-linearizability condition holds true, for instance, if between the eigenvalues $\alf_{1,2}$ there are no resonant relations $\alf_i = n_1\alf_1 + n_2\alf_2$ with integers $n_{1,2} \geq 0$, $n_1+n_2 \ge 2$
(Sternberg–Chen Theorem, see e.g., \cite{AI, Hartman}).
The proof presented in \cite{Dav85} is done in $C^{\infty}$ category, but is valid in $C^{\om}$ as well
(the requirement of $C^{\infty}$-linearizablity should be replaced with $C^{\om}$-linearizablity),
see also the recent paper \cite{Boga}.


\section{Singular points of the geodesic flow}

In addition to the isotropic surface $\FF$, the vector field $\vec V$ given by \eqref{3} has one more evident
invariant surface~-- the {\it vertical} surface
$$
\ov S_0 = \{ (q,p) \ \ q=(x,y) \in S_0, \ \ p \in \PR\}.
$$
The restriction of the field \eqref{3} to $\ov S_0$ is vertical at almost all points
(except for the points where $M=0$, and the field vanishes). Hence the surface $\ov S_0$ is filled with
vertical integral curves of the filed \eqref{3} and its singular points.

Singular points of the field \eqref{3} are given by two equations:
\begin{equation}
\Del(x,y) = 0 \  \ \ \textrm{and} \ \ \  M(x,y,p)=0,
\label{4}
\end{equation}
and consequently, they are not isolated, but form a curve (or curves) in $PTS$.
Algebraically, this property can be expresses in the following form:
all components of the vector field \eqref{3} belong to the ideal $I$
(in the ring of smooth functions) generated by two of them,
namely, $I =\<\Del, M\>$.

\begin{remark}
{\rm
The fact that the {\it horizontal} generator $\Del(x,y)$ of the field \eqref{3}
does not depend on $p$ and the {\it vertical} generator $M(x,y,p)$ is a cubic polynomial in $p$,
plays a crucial role in a general geometrical context, e.g.,
in the framework of Cartan's theory of the projective connection \cite{Aminova, A-Geom}.
}
\label{Rem1}
\end{remark}

Let us list those of the properties of the field \eqref{3} that we are going to use:
\begin{itemize}
\item
Singular points of the field \eqref{3} are given by equation \eqref{4} and
form a curve (or several curves) in $PTS$.
\item
The spectrum of the linearization of the field \eqref{3} at every singular point contains one zero eigenvalue
and two real eigenvalues $\la_{1,2}$, which vanish (simultaneously) at those points where the cubic polynomial $M(q,p)$ has a double root $p$. The latter condition is equivalent to the direction $p$ is tangent
to $S_0$ at the point $q$.
\item
For every point $q \in S_0$ and any $p \in \PR$ such that $M(q,p) \neq 0$ there exists
a unique integral curve of the field \eqref{3} that passes through the point $(q,p)$~--
a vertical straight line, whose projection on $S$ is not a geodesic.
Consequently, the vertical surface $\ov S_0$ is an invariant surface of \eqref{3}.
\item
Geodesics cannot enter a point $q \in S_0$ in arbitrary tangential direction, but only in {\it admissible}
directions $p$ that satisfy the condition $M(q,p) = 0$.
\item
The isotropic direction $p_0(q)$ given by formula \eqref{40} is admissible at every point $q \in S_0$, i.e.,
$M(q,p_0(q)) = 0$ for all $q \in S_0$.
\end{itemize}

\begin{figure}
\begin{center}
\includegraphics[height=395pt]{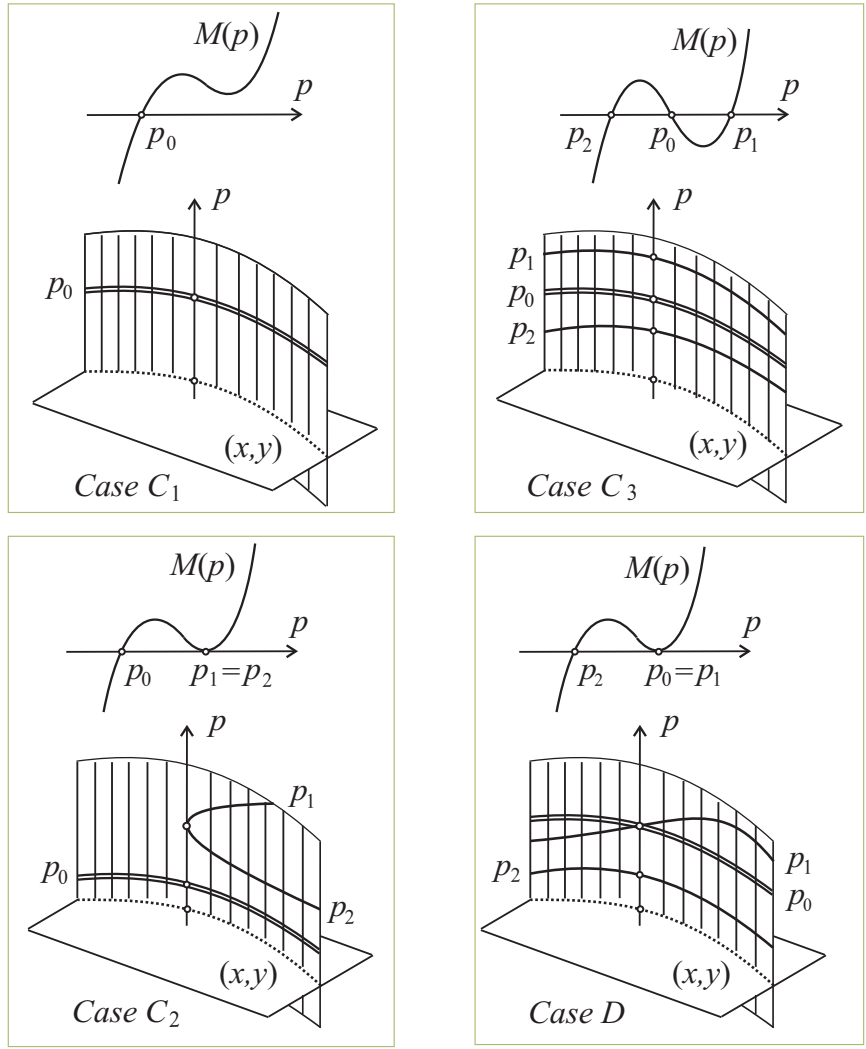}
\caption{
Real roots of the cubic polynomial $M(p)$ and the set of singular point of the field \eqref{3}.
The double line presents $\{ (q,p_0), \ q \in S_0\}$, the bold lines present $\{(q,p_i), \ q \in S_0\}$, $i=1,2$,
the dotted line presents~$S_0$.
}
\label{fig4}
\end{center}
\end{figure}

Depending on the roots of the cubic polynomial $M$ (see Fig.~\ref{fig4}), we have four cases:
\begin{itemize}
\item
$C_1$ : the isotropic direction $p_0$ is a unique real root of $M$,
\item
$C_2$ : $M$ has a simple root $p_0$ and a double non-isotropic real root $p_{1}=p_{2}$,
\item
$C_3$ : $M$ has three simple real roots: isotropic $p_0$ and non-isotropic $p_{1}, p_{2}$,
\item
$D \phantom{.}$ : the isotropic double root $p_0=p_1$ and a simple non-isotropic root $p_{2}$.
\end{itemize}

If $\Real \la_{1,2} \neq 0$, the set $W$ is the center manifold of the field,
and the restriction of the field to $W$ is identically zero.
Hence in a neighborhood of every singular point where $\Real \la_{1,2} \neq 0$, the phase portrait
of the field has a very simple topological structure.
Indeed, the reduction principle \cite{AI, HPS} asserts that
the germ of the field is orbitally topologically equivalent
to the direct product of the standard 2-dimensional node (if $\Real \la_{1,2}$ have the same sign)
or saddle (if $\Real \la_{1,2}$ have different signs) and 1-dimensional zero vector field.
However, the topological classification is not enough.

The paper \cite{Rem06} presents finite-smooth local normal forms of such fields,
\cite{Rem-Tari} contains a brief survey (Appendix A) on the smooth and $C^{\omega}$ classifications.
These results allow to establish smooth local normal forms of the field \eqref{3} at all singular points
$(q,p_i)$, $q \in S_0$, where $p_i$ is a simple real root of $M(q,p)$.
This gives the description of geodesics that enter a degenerate point with all possible admissible directions for the cases $C_1$, $C_3$. To study geodesics with the isotropic admissible direction in the cases $C_2$ and $D$,
one can use a blow-up procedure.

Choosing appropriate local coordinates, we shall further assume that in a neighborhood of the point $q \in S_0$,
equation \eqref{2} has the form $dx^2 = y \, dy^2$ in the case $C$ and \eqref{DD} in the case $D$.
Consequently, the discriminant curve $S_0$ is the axis $y=0$ in the case $C$ and the parabola
$y - \eps x^2 = 0$ in the case $D$.
Since multiplication the metric by the factor~$-1$ does not change the geodesic flow,
without loss of generality, assume that $y > \eps x^2$ and $y < \eps x^2$ (including the case $\eps=0$)
are Lorentzian and Riemannian domains, respectively.

From now on, we shall consider geodesics outgoing from a degenerate point $q \in S_0$ with the isotropic admissible
direction $p_0(q)$ as semitrajectories starting from $q$.
We distinguish geodesics outgoing into the Lorentzian (resp. Riemannian) domains using the superscript $+$ (resp. $-$).
Let us clarify this with the following example.

\begin{exam}
{\rm
For the metric $ds^2 = dx^2 - y dy^2$, the discriminant curve $S_0 = \{ y=0\}$
divides the plane into the Lorentzian ($y>0$) and Riemannian ($y<0$) domains.
Formula \eqref{mu} yields $M(q,p)=p^2$, and we have the case $C_2$.

At every degenerate point $q \in S_0$ there exist two
admissible directions: $p_1 = 0$ (non-isotropic, double root) and $p_0=\infty$ (isotropic).
To see that the direction $p_0=\infty$ is admissible, it is convenient to interchange $x$ and $y$.
In the new coordinates $\ox=y$, $\oy=x$, $\op=1/p$,
the polynomial $M(q,\op)=-\op$ has the root $\op=0$.

The corresponding field \eqref{3} has a unique integral curve $y = 0$ that pass through every point $q \in S_0$
with tangential direction $p_1=0$. Substituting $y=0$ directly in \eqref{ELA}, one can see that $y=0$
is an extremal of the action functional and its natural parametrization is given by the equation $\ddot x = 0$.
Moreover, given degenerate point $q \in S_0$ there exists a one-parameter family of geodesics
outgoing from $q$ with the tangential direction $p=\infty$.
For instance, consider the family $\Gamo$ of geodesics $\gama$, $\alf \in \bR$,
outgoing from the origin; see Fig.~\ref{fig7} (right). They can be presented in the agreed upon way as
\begin{equation}
\gama = \left \{ \
\begin{aligned}
\gamap : \ &x = \alf y^{\frac{3}{2}}, \ \ \ \ \ \, y \ge 0, \\
\gamam : \ &x = \alf (-y)^{\frac{3}{2}}, \ y \le 0. \\
\end{aligned}
\right.
\label{Gamma}
\end{equation}
}
\label{Ex1}
\end{exam}


\subsection{The case $C$}

The linearization of the field \eqref{3} at every singular point $(q,p_i)$, $q \in S_0$, $i=0,1,2$,
has the spectrum $(\la_1,\la_2,0)$ with non-zero real eigenvalues $\la_{1,2}$.
Moreover, at a singular point $(q,p_0)$ corresponding to the isotropic admissible direction
the resonant relation $\la_1 = 2\la_2$ holds.
On the other hand, at a singular point $(q,p_i)$, $i=1,2$, corresponding to non-isotropic admissible direction
the resonant relation $\la_1 + \la_2 = 0$ holds.\footnote{
The relation $\la_1 = 2\la_2$ is a corollary of $\la_1 + \la_2 + \la_3 = r\lambda_1$
with $r=\tfrac{3}{2}$ and $\la_3=0$, see Theorem~\ref{T0} and formula \eqref{5-Oct-2016}.
The relation $\la_1 + \la_2 = 0$ follows form the fact that the field $\vec W$ is divergence-free
and the function $F$ does not vanish in a neighborhood of $(q,p_i)$, $i=1,2$.
}
Using the smooth classification of vector fields with non-isolated singular points
(see e.g. \cite{Rem-Tari}, Appendix~A), we have the following results.

The germ of the field \eqref{3} at any point $(q,p_0)$, $q \in S_0$, has $C^\infty$ orbital normal form
\begin{equation}
2\xi \frac{\pa}{\pa \xi} + \eta \frac{\pa}{\pa \eta}  + 0 \frac{\pa}{\pa \zeta}
\label{401}
\end{equation}
with the first integrals $I_1=\xi/\eta^2$ and $I_2=\zeta$.
The germ of the field \eqref{3} at any point $(q,p_i)$, $q \in S_0$, $i=1,2$,
has $C^\infty$ orbital normal form
\begin{equation}
\xi \frac{\pa}{\pa \xi} - \eta \frac{\pa}{\pa \eta}  + \xi \eta \frac{\pa}{\pa \zeta}
\label{403}
\end{equation}
with the first integral $I= \xi \eta $.
One can see that to every singular point of the field \eqref{401} corresponds a one-parameter family
of integral curves passing through this point, while to every singular point of the field \eqref{403}
correspond only two integral curves. Projecting the integral curves down, we obtain the following results.

\begin{figure}
\begin{center}
\includegraphics[width=370pt,height=210pt]{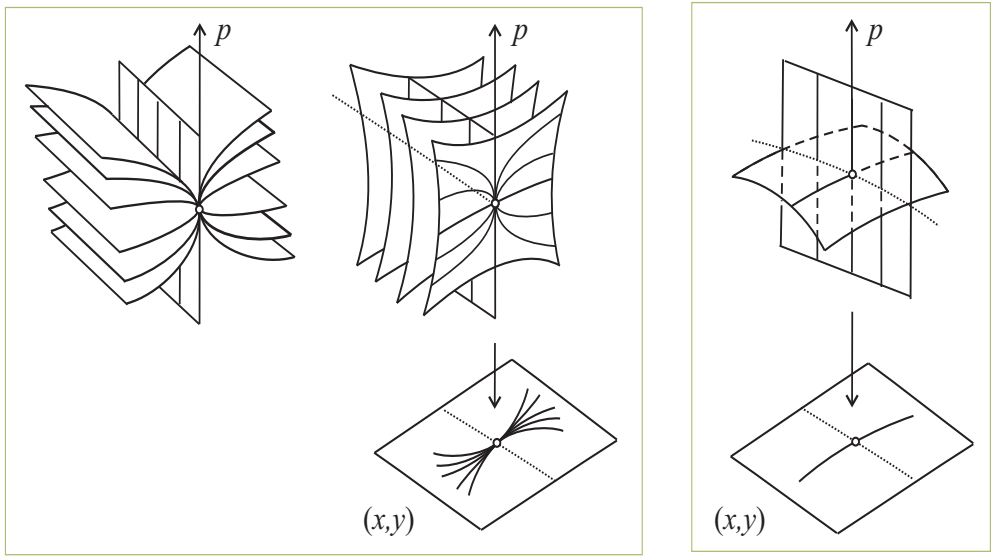}
\vspace{-1ex}
\caption{
The cases $C_1$, $C_3$.
Left panel:
two invariant foliations of the field \eqref{3} near the point $(q,p_0)$, $q \in S_0$.
Right panel:
the invariant leaves of the field \eqref{3} passing through the point $(q,p_i)$, $q \in S_0$, $i=1,2$.
The dotted lines present the set of singular points of the field \eqref{3} and
its projection, the discriminant curve.
}
\label{fig6}
\end{center}
\end{figure}

\begin{theorem}[\cite{Rem-Pseudo, Rem15}]
\label{T1}
Suppose that $C$ holds true.
Then to the isotropic direction $p_0$ corresponds a one-parameter family $\Gamo$ of geodesics
outgoing from the point $q$.
There exist smooth local coordinates centered at $q$ such that
the discriminant curve $S_0$ coincides with the $x$-axis,
the isotropic direction $p_0(q)=\infty$
and the geodesics $\gamapm \in \Gamo$ are semi-cubic parabolas
\begin{equation}
x= \alf \tau^3 X_{\alf}^{\pm}(\tau), \quad
y= \tau^2 Y_{\alf}^{\pm}(\tau), \quad \alf \ge 0,
\label{402}
\end{equation}
where $X_{\alf}^{\pm}, Y_{\alf}^{\pm}$ are smooth functions, $X_{\alf}^{\pm}(0)=1$, $Y_{\alf}^{\pm}(0)=\pm 1$.
\end{theorem}

\begin{theorem}[\cite{Rem-Pseudo, Rem15}]
\label{T2}
Suppose that $C_3$ holds true.
Then to each admissible direction $p_i$, $i=1,2$,
corresponds a unique geodesic passing through the point $q$.
Both these geodesics are smooth and timelike.
\end{theorem}

In the left panel of Fig.~\ref{fig6} we present the invariant foliations of the field \eqref{3}
in a neighborhood of the point $(q,p_0)$, $q \in S_0$, that correspond to the first integrals
$I_1=\xi/\eta^2$ (left) and $I_2=\zeta$ (right) of the normal form \eqref{401}.
Intersection of these foliations gives the family of integral curves of \eqref{3}.
The family $\Gamo$ of the geodesics \eqref{402} is obtained (by the projection $PTS \to S$) from
the family of integral curves of the field \eqref{3} that pass through its singular point $(q,p_0)$.
The subfamily $\Gamo^+ \subset \Gamo$ of the geodesics \eqref{402} outgoing into the Lorentzian semiplane,
contains timelike, spacelike, and isotropic geodesics.

In the right panel of Fig.~\ref{fig6} we present
those of the leaves of the invariant foliation of the field \eqref{3}
in a neighborhood of the point $(q,p_i)$, $q \in S_0$, $i=1,2$, that pass through $(q,p_i)$.
This foliation corresponds to the first integral $I=\xi\eta$ in the normal form \eqref{403},
and the leaves passing through $(q,p_i)$ coincide with the planes $\xi=0$ and $\eta=0$,
while none of the remaining leaves contains singular points of \eqref{3}.
One of these leaves coincides with the vertical surface $\ov S_0$ filled with vertical
integral curves whose projection on $S$ are points of $S_0$.
Another invariant surface is filled with non-vertical integral curves,
through every point $(q,p_i)$, $q \in S_0$, there pass exactly one curve.

\begin{figure}
\begin{center}
\includegraphics[width=350pt,height=130pt]{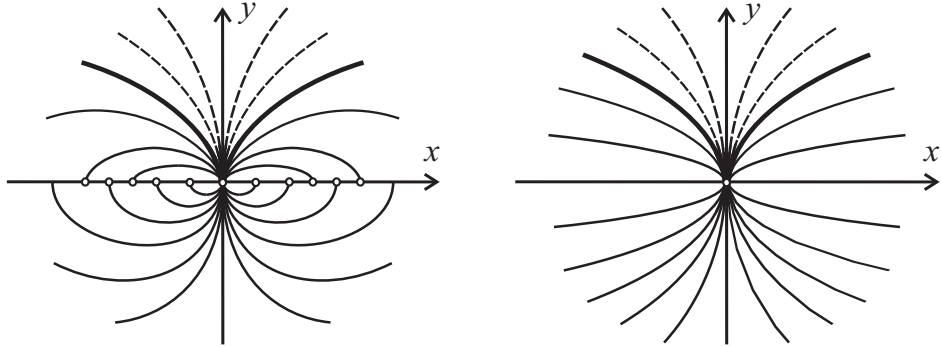}
\vspace{-1ex}
\caption{
Two examples of the family $\Gamo$ in Theorem~\ref{T1}.
Geodesics in the metrics $e^y dx^2 - y dy^2$ (left) and $dx^2 - y dy^2$ (right)
outgoing from $q = 0$.
Timelike, spacelike, and isotropic geodesics are depicted as solid, dashed, and bold solid lines respectively.
}
\label{fig7}
\end{center}
\end{figure}

\begin{figure}
\begin{center}
\includegraphics[width=370pt,height=125pt]{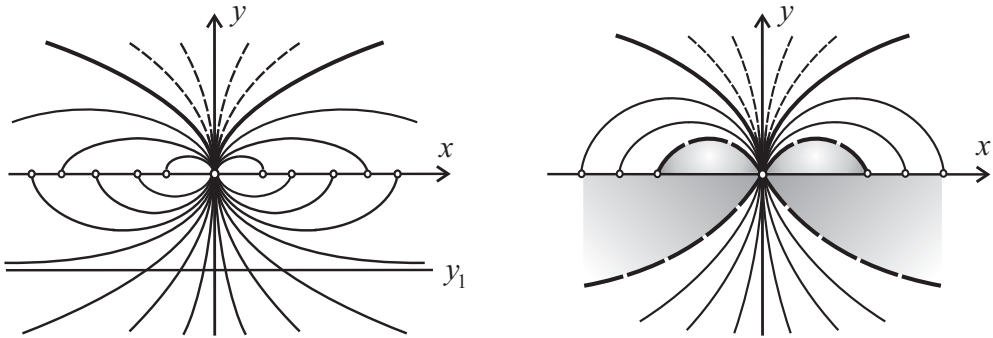}
\vspace{-1ex}
\caption{
Geodesics in the metrics $ds^2 = a(y)dx^2 - ydy^2$ with $a(y)=1+(y + y_1)^2$ outgoing from $q=0$.
The case $C_1$ ($y_1<0$, left) and the case $C_3$ ($y_1>0$, right).
Timelike, spacelike, and isotropic geodesics are depicted as solid, dashed, and bold solid lines respectively.
Geodesics passing through $0$ with non-isotropic admissible directions (right) are depicted as long-dashed bold lines.
The grey domains do not contain geodesics passing though~$0$.
}
\label{fig9}
\end{center}
\end{figure}

\begin{exam}
{\rm
To illustrate the above, return to Example~\ref{Ex1}.
In the coordinates $\ox=y$, $\oy=x$, $\op=1/p$,
the equation of isotropic lines coincides with Cibrario normal form.
After multiplication by $-1$, the corresponding vector field \eqref{3} reads
\begin{equation}
\vec V = 2\ox \biggl(\frac{\pa}{\pa \ox} + \op \frac{\pa}{\pa \oy}\biggr) + \op \frac{\pa}{\pa \op}.
\label{9-Oct-2016}
\end{equation}

It is easy to check that the field \eqref{9-Oct-2016} possesses the invariant foliation $\ox = c \op^2$,
which includes, in particular, the vertical surface $\ov S_0$ (for $c=0$),
the isotropic surface (for $c=1$).
This foliation is presented in the left side of the left panel of Fig.~\ref{fig6}.
The restriction of the field \eqref{9-Oct-2016} to every invariant leaf $\ox = c \op^2$ reads
$2c\op^3 \frac{\pa}{\pa \oy} + \op \frac{\pa}{\pa \op}$.
Canceling the factor $\op$, we obtain the non-singular field
$2c\op^2 \frac{\pa}{\pa \oy} + \frac{\pa}{\pa \op}$, whose integral curves are presented in Fig.~\ref{fig5} (left).
Fixing a degenerate point $q \in S_0$, in going through all invariant leaves $\ox = c \op^2$ and projecting down,
we obtain the family \eqref{Gamma} of geodesics $\gamap$ (for $c>0$) and $\gamam$ (for $c<0$)
presented in Fig.~\ref{fig7} (right).\footnote{
The attentive reader may remark that this invariant foliation contains also the leaf $\op=0$,
which can be considered as the limiting case for $c \to \infty$.
The restriction of \eqref{9-Oct-2016} to this leaf is filled with integral curves parallel to the $\ox$-axis.
This gives the family of geodesics $x = \const$, which are
the limiting case of the semi-cubic parabolas \eqref{402}: the two branches are glued together.
}
}
\label{Ex2}
\end{exam}

\begin{remark}
{\rm
If the pseudo-Riemannian metric on the surface $S$ is induced by
the pseudo-Euclidean metric $dX^2+dY^2-dZ^2$ of the ambient space (see the example above),
the difference between the cases $C_1$ and $C_3$ has a graphical interpretation.
Namely, $C_1$ and $C_3$ correspond to positive and negative Gaussian curvature of the surface $S$
calculated in the Euclidean metric $dX^2+dY^2+dZ^2$.
}
\label{Rem3}
\end{remark}

\begin{figure}
\begin{center}
\includegraphics[height=138pt]{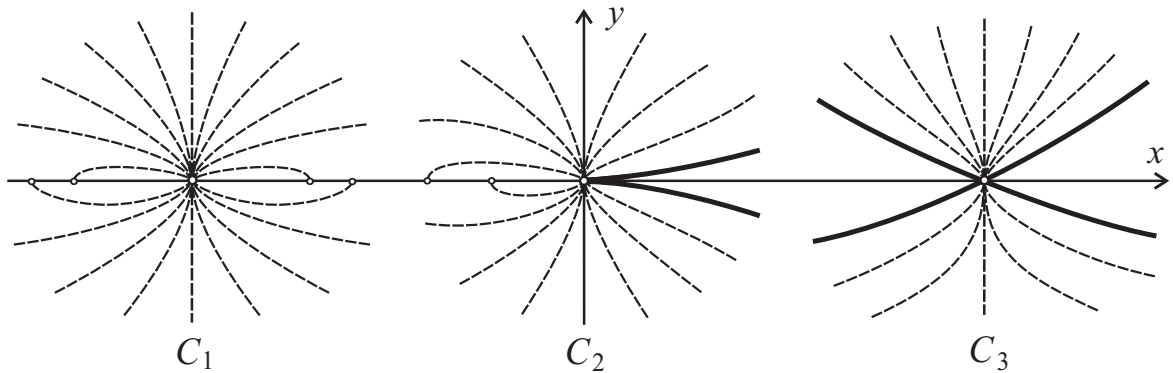}
\caption{
Geodesics in the metrics \eqref{29-11-2016} outgoing from three different degenerate points:
$C_1$ ($x<0$, left), $C_2$ ($x=0$, center), and $C_3$ ($x>0$, right).
Here timelike, spacelike, and isotropic geodesics outgoing from $q \in S_0$ with the isotropic direction
$p_0=\infty$ are depicted as dashed lines, while the geodesics outgoing from $q \in S_0$
with non-isotropic admissible direction $p_1=0$ are depicted as bold lines.
}
\label{figC2}
\end{center}
\end{figure}

\begin{theorem}
\label{C2}
Suppose that $C_2$ holds true.
Generically, the point $q$ locally separates the curve $S_0$ in two parts,
filled with $C_1$ and $C_3$ points, respectively, and there exist smooth local coordinates centered at $q$
such that the metric has the form
\begin{equation}
ds^2 = a(x,y)\, dx^2 + ye(x,y)\, dy^2, \ \ a(0) \neq 0, \ e(0) \neq 0, \ a_y(0)=0, \ a_{xy}(0) \neq 0.
\label{29-11-2016}
\end{equation}

Then to the double admissible direction $p_1=p_2$ corresponds a unique geodesic passing through the point $q$,
a semicubic parabola with branches outgoing from $q$ into the Lorentzian and Riemannian domains
(depicted as long-dashed line in Fig.~\ref{figC2}, center).
\end{theorem}

The proof is not published yet.
In Example~\ref{Ex1} considered above, we deal with a non-generic case $C_2$, since the condition
$a_{xy}(0) \neq 0$ in \eqref{29-11-2016} does not hold true.
This leads to the geodesic $y=0$ instead of a semicubic parabola mentioned in Theorem~\ref{C2}.

\subsection{The case $D$}

The cubic polynomial $M$ at $q \in S_0$
has the isotropic double root $p_0 = p_1$ and a simple non-isotropic root $p_2$.
For the admissible direction $p_2$, the analogous assertion to Theorem~\ref{T2} holds true:
the germ of the field \eqref{3} at $(q,p_2)$ has $C^\infty$ normal form \eqref{403}, and
to the direction $p_2$ corresponds a unique smooth geodesic passing through the point $q$.
However, the study of geodesics with the isotropic direction is more complicated.

A special feature of the case $D$ is that the linear part of the germ \eqref{3}
at $(q,p_0)$, $q \in S_0$, has three zero eigenvalues.
This prevents the possibility to obtain a normal form
similar to \eqref{401} in Theorem~\ref{T1} or similar to \eqref{403} in Theorem~\ref{T2}.
Moreover, in this case even the reduction principle does not allow to establish the topological
normal form of this filed, since the {\it center subspace}\footnote{
The center subspace $T_c$ of a vector filed $\vec V$ at its singular point $0$
is spanned by the generalized eigenvectors of the linearization of $\vec V$ at $0$
corresponding to the eigenvalues $\la$ with $\Real \la = 0$.
}
of the germ \eqref{3} at $(q,p_0)$ coincides with the whole tangent space, see \cite{AI}.
However, using appropriate blowing up procedure, one can reduce the germ \eqref{3} at $(q,p_0)$
to a smooth vector field with non-zero spectrum and study the obtained vector field using the standard methods.

Further we always assume that in the cases $D_s$ and $D_n$ the following genericity condition holds true:
there are no non-trivial integer relations
\begin{equation*}
n_1\alf_1 + n_2\alf_2 + n_3 \alf_3 = \alf_j, \ \ \,
n_1 + n_2 + n_3 \ge 1, \ \ n_i \in \bZ_+, \ \ \, j=1,2,3,
\end{equation*}
where $\alf_{1,2}$ are the eigenvalues of the linearization of the vector field \eqref{6-Oct-2016} at $(q,p_0)$
and $\alf_3=2$. This condition implies the germ of a vector field obtained from \eqref{3}
by the blowing up procedure is linearizable, as well as the germ of the field \eqref{6-Oct-2016}.

\subsubsection{The cases $D_n$ and $D_f$}

\begin{figure}
\begin{center}
\includegraphics[height=5.7cm]{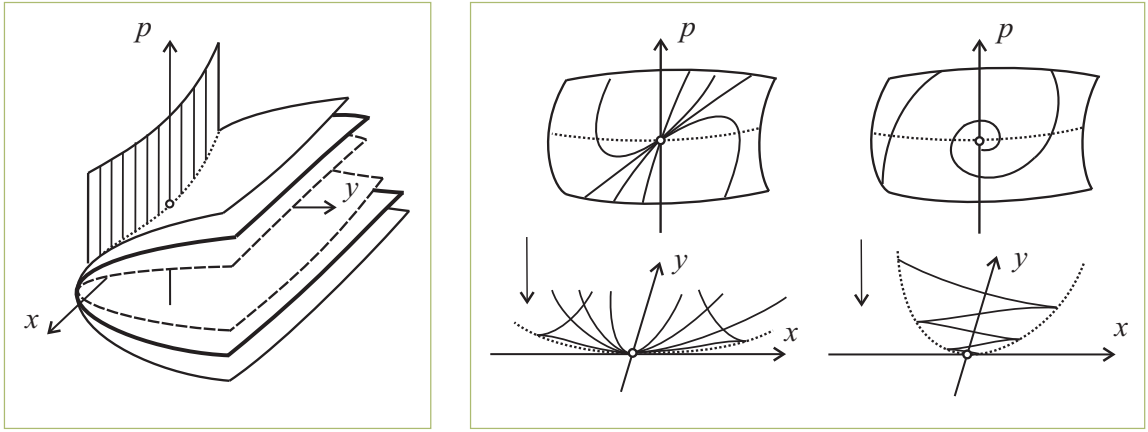}
\end{center}
\vspace{-1ex}
\caption{
The cases $D_n$ and $D_f$.
The left panel: the invariant foliation $\{\FF_{\alf}\}$ of the field \eqref{3} above the Lorentzian domain
(the isotropic surface $\FF = \FF_0$ is depicted as bold). Here the leaves filled with
timelike, spacelike and isotropic geodesics are depicted as solid, dashed and bold solid lines, respectively.
The right panel: integral curves of the restriction of field \eqref{3} to $\FF_{\alf}$ and their projections
(the case $D_n$ on the left and $D_f$ on the right).
The criminant and the discriminant curve are depicted as dotted lines.
}
\label{pic7}
\end{figure}

In a neighborhood of the considered point $(q,p_0)$, $q \in S_0$,
the field \eqref{3} above the Lorentzian domain
has an invariant foliation $\{\FF_{\alf}\}$ presented in the left panel of Fig.~\ref{pic7}.
Here the invariant leaf $\FF_{0}$ coincides with the isotropic surface $\FF$.
The invariant leaves above the Riemannian domain are not depicted, since
they contain no integral curves that pass through $(q,p_0)$.

The linear part of the restriction of the field \eqref{3} to every invariant leaf $\FF_{\alf}$
at its singular point $(q,p_0)$
is equal to \eqref{6-Oct-2016} multiplied by a smooth scalar function $\si_{\alf}$
vanishing along the criminant. Therefore, the restriction of the field \eqref{3}
to every invariant leaf $\FF_{\alf}$ has the local phase portrait of the same type: node or focus.
See Fig.~\ref{pic7} (right panel).
In going through all invariant leaves $\FF_{\alf}$ and projecting the integral curves down,
we obtain the following result.

\begin{theorem}[\cite{Rem-Tari}]
\label{T3}
Let the case $D_n$ or $D_f$ holds true.
Then to the isotropic direction $p_0$ corresponds a two-parameter family $\Gamo$ of
$C^2$-smooth geodesics $\gamap$ outgoing from $q$ into the Lorentzian domain, while
there are no geodesics outgoing from $q$ into the Riemannian domain.
Given $\alf$, the geodesics $\gamabp \in \Gamo$ with fixed $\alf$ and varying $\bet$
are projections of the integral curves from the leaf $\FF_{\alf}$;
see Fig.~\ref{pic7}, center for $D_n$ and right for $D_f$.
The geodesics $\gamabp \in \Gamo$ are timelike if $\alf<0$, spacelike if $\alf>0$ and isotropic if $\alf=0$.
\end{theorem}

\subsubsection{The case $D_s$}

\begin{figure}
\begin{center}
\includegraphics[height=5.3cm]{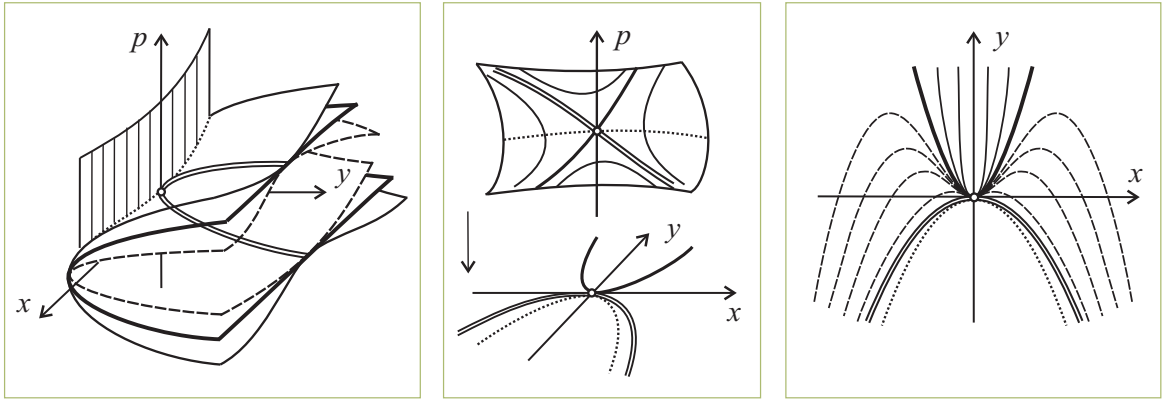}
\end{center}
\vspace{-1ex}
\caption{
The case $D_s$.
On the left: invariant foliation $\{\FF_{\alf}\}$ of the field \eqref{3} above the Lorentzian domain
(the isotropic surface $\FF = \FF_0$ is depicted as bold).
Center: integral curves of the field \eqref{3} on an invariant leaf $\FF_{\alf}$.
On the right: geodesics outgoing from the point $q \in S_0$.
Timelike and spacelike geodesics are depicted as solid and dashed lines, respectively.
The bold solid and the double solid lines present two isotropic geodesics.
The criminant and the discriminant curve are depicted as dotted lines.
}
\label{pic6}
\end{figure}

In a neighborhood of the considered point $(q,p_0)$, $q \in S_0$,
the field \eqref{3} above the Lorentzian domain
has an invariant foliation $\{\FF_{\alf}\}$ presented in the left panel of Fig.~\ref{pic6}.
Here the invariant leaf $\FF_{0}$ coincides with the isotropic surface $\FF$.
The invariant leaves above the Riemannian domain are not depicted, since
they contain no integral curves that pass through $(q,p_0)$.

The linear part of the restriction of the field \eqref{3} to every invariant leaf $\FF_{\alf}$
at its singular point $(q,p_0)$
is equal to \eqref{6-Oct-2016} multiplied by a smooth scalar function $\si_{\alf}$
vanishing along the criminant. Therefore, the restriction of the field \eqref{3}
to every invariant leaf $\FF_{\alf}$ has a saddle at $(q,p_0)$.
See Fig.~\ref{pic6} (right panel).

\begin{theorem}[\cite{Rem-Tari}]
\label{T4}
Let the case $D_s$ holds true.
Then to the isotropic direction $p_0$ corresponds a one-parameter family $\Gamo$ of
$C^2$-smooth geodesics outgoing from $q$ into the Lorentzian domain, while
there are no geodesics outgoing from $q$ into the Riemannian domain.
There exist smooth local coordinates centered at $q$ such that
$S_0$ is the parabola $y = \eps x^2$ and
the geodesics $\gamap \in \Gamo$ outgoing from $q$ have the form
\begin{equation}
y = \frac{\eps_1}{2}x^2 + Y_{\alf}(x), \ \ \,  Y_{\alf}(x)=o(x^2), \, \ \  \alf \in \mathbb R,
\label{26}
\end{equation}
together with one additional isotropic geodesic
\begin{equation}
y = \frac{\eps_2}{2}x^2 + Y(x), \ \ \,  Y(x)=o(x^2),
\label{27}
\end{equation}
where $\eps_1 \eps_2 = \eps$, $\eps_1 + \eps_2 = \frac{1}{2}$, $\eps_1 > \tfrac{1}{2}$, $\eps_2 < 0$.
Geodesics \eqref{26} are timelike if $\alf<0$, spacelike if $\alf>0$, isotropic if $\alf=0$;
see Fig.~\ref{pic6}, right.
\end{theorem}

\begin{figure}
\begin{center}
\includegraphics[height=5.4cm]{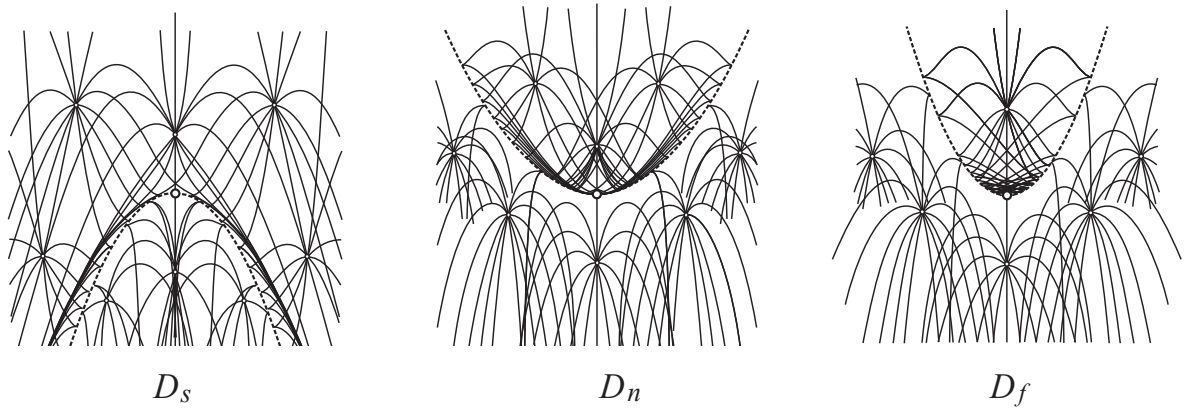}
\end{center}
\vspace{-1ex}
\caption{
Computer generated geodesics (solid lines) of the metric $dy^2 + (\eps x^2 - y) dx^2$
with $\eps < 0$ (the case $D_s$), $0 < \eps < \tfrac{1}{16}$ (the case $D_n$), and
$\eps > \tfrac{1}{16}$ (the case $D_f$).
The parabola depicted as dotted line is the discriminant curve $S_0$.
}
\label{pic8}
\end{figure}

It is interesting to note that invariant foliations in the cases $D_n$, $D_f$ and $D_s$ have
the different topological structures (compare the left panels of Figures~\ref{pic7} and~\ref{pic6}).
In the cases $D_n$, $D_f$ all invariant leaves intersect on the criminant only,
while in the case $D_s$ they intersect on the criminant (dotted line) and on the double line,
whose projection is the isotropic geodesic \eqref{27}.


\subsection{Example: Clairaut type}

It is of interest to observe an important difference between the families $\Gamo$ in the cases $C_1$, $C_3$ and $D$.
In the cases $C_1$, $C_3$, the family $\Gamo$ is symmetric with respect to $S_0$ in the following sense:
it contains an infinite number of geodesics $\gamap \in \Gamo$ outgoing into the Lorentzian domain
and an infinite number of geodesics $\gamam \in \Gamo$ outgoing into the Riemannian domain.
On the contrary, in the case $D$, the family $\Gamo$ is non-symmetric: it contains
an infinite number of geodesics $\gamap \in \Gamo$ outgoing into the Lorentzian domain
and no geodesics $\gamam \in \Gamo$ outgoing into the Riemannian domain.

To understand this phenomenon better, consider the case when the isotropic direction $p_0$
is tangent to the curve $S_0$ {\it at all points} $q \in S_0$, for instance,
the metric $dy^2 + (\eps x^2 - y) dx^2$.
The equation of geodesics in the metric $ds^2 =  dy^2 - y dx^2$ can be studied using qualitative methods,
see \cite{Rem15} (Section~3).
The Lagrangian of the length functional $L=\sqrt{p^2-y}$ does not depend on the variable $x$,
hence the field \eqref{3} possesses the energy integral $H = L-pL_p$.
After evident transformations, equation $H = \const$ can be reduced to
\begin{equation}
p^2 = y - \alf y^2, \ \ \, \alf \in \bR,
\label{43000}
\end{equation}
which is a family of implicit differential equations {\it of Clairaut type} \cite{Dav-Japan}.

Every (unparametrized) geodesic in the metric $ds^2 =  dy^2 - y dx^2$ is a solution of equation \eqref{43000}.
Conversely, every solution of \eqref{43000} is a geodesic except the horizontal lines $y \equiv \const$,
each of which is the envelop of the family of integral curves of \eqref{43000} for a given $\alf$
(see \cite{Rem15}).
For instance, the value $\alf=0$ corresponds to the isotropic surface $p^2 = y$ (a parabolic cylinder)
and gives, in particular, the isotropic geodesic $y=\frac{1}{4}x^2$ passing through the origin.

\begin{figure}
\begin{center}
\includegraphics[height=4.6cm]{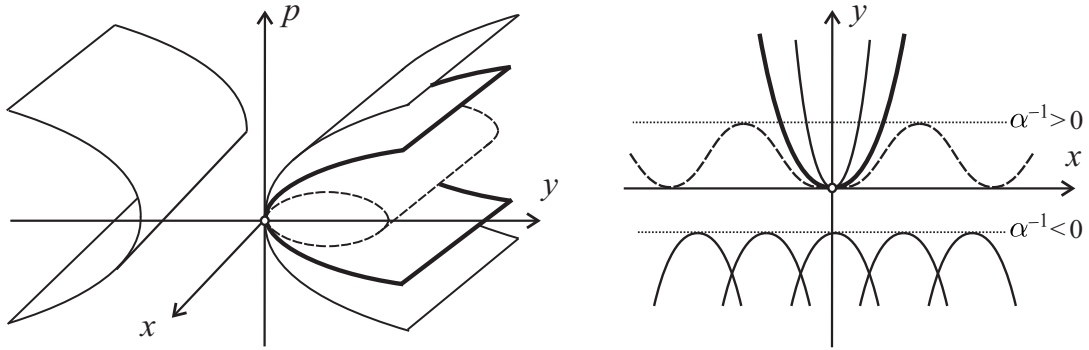}
\end{center}
\vspace{-1ex}
\caption{
The invariant foliation $p^2 = y-\alf y^2$ in the $(x,y,p)$-space (left)
and the corresponding geodesics (right).
Timelike, spacelike and isotropic geodesics are solid, dashed and bold solid lines, respectively.
}
\label{pic5}
\end{figure}

For determining non-isotropic geodesics, observe that every invariant surface \eqref{43000}
is a cylinder whose generatrices are parallel to the $x$-axis
and the base is an ellipse (if $\alf>0$) or a hyperbola (if $\alf<0$).
In the latter case, the hyperbolic cylinder $p^2 = y-\alf y^2$ consists of two connected components:
{\it positive} and {\it negative} lying in the domains $y \ge 0$ and $y \le \alf^{-1}$, respectively.
Positive components of the hyperbolic cylinders ($\alf<0$) together with all other cylinders ($\alf \ge 0$)
form an invariant foliation over the Lorentzian domain $y>0$.
Negative components of the hyperbolic cylinders
form an invariant foliation over the Riemannian domain $y<0$;
they do not intersect the plane $y=0$, and consequently,
do not contain integral curves whose projections to the $(x,y)$-plane are geodesics passing through the $x$-axis.
See Fig.~\ref{pic5} (left).

Thus to every $\alf \ge 0$ corresponds a geodesic $\gamap \in \Gamo$ which is timelike if $\alf>0$
or isotropic if $\alf=0$.
To every $\alf < 0$ corresponds a spacelike geodesic $\gamap \in \Gamo$, whose lift belongs
to the positive component of the hyperbolic cylinder $p^2 = y-\alf y^2$.
In contrast to this, the negative component of the same cylinder is filled with integral curves
of the field \eqref{3} whose projections on the $(x,y)$-plane are separated from the $x$-axis by
the horizontal strip $\alf^{-1}<y<0$.
Therefore, there are no geodesics outgoing into the Riemannian domain.
See Fig.~\ref{pic5}, right.


\medskip

{\bf Acknowledgement}. \
The publication was supported by the Russian Foundation for Basic Research
(research projects 16-01-00766, 17-01-00849).


\end{document}